# *P*-VARIATION OF STRONG MARKOV PROCESSES

By Martynas Manstavičius

*University of Connecticut*

Let $\xi_t$, $t \in [0,T]$, be a strong Markov process with values in a complete separable metric space $(X, \rho)$ and with transition probability function $P_{s,t}(x, dy)$, $0 \leq s \leq t \leq T$, $x \in X$. For any $h \in [0, T]$ and $a > 0$, consider the function $\alpha(h, a) = \sup\{P_{s,t}(x, \{y : \rho(x,y) \geq a\}) : x \in X, 0 \leq s \leq t \leq (s+h) \wedge T\}$. It is shown that a certain growth condition on $\alpha(h,a)$, as $a \downarrow 0$ and $h$ stays fixed, implies the almost sure boundedness of the *p*-variation of $\xi_t$, where $p$ depends on the rate of growth.

**1. Introduction.** Let $\xi_t$, $t \in [0,T]$, be a strong Markov process defined on some complete probability space $(\Omega, \mathcal{F}, P)$ and with values in a complete separable metric space $(X, \rho)$. Denote the transition probability function of $\xi_t$ by $P_{s,t}(x, dy)$, $0 \leq s \leq t \leq T$, $x \in X$. For any $h \in [0, T]$ and $a > 0$, consider the function

$$\alpha(h, a) = \sup\{P_{s,t}(x, \{y : \rho(x,y) \geq a\}) : x \in X, 0 \leq s \leq t \leq (s+h) \wedge T\}.$$

The behavior of $\alpha(h, a)$ as a function of $h$ gives sufficient conditions for regularity properties of the trajectories of the process $\xi_t$. As Kinney (1953) showed, $\xi_t$ has an almost surely càdlàg version if $\alpha(h, a) \to 0$ as $h \to 0$ for any fixed $a > 0$, and an almost surely continuous version if $\alpha(h, a) = o(h)$ as $h \to 0$ for any fixed $a > 0$. Dynkin (1952) obtained the same result but under slightly stronger conditions. The main goal of this paper is to establish a connection between the arguments of $\alpha(h,a)$ and variational properties of paths of $\xi_t$, in particular, *p*-variation which in many instances is very useful, as discussed in Dudley (1992), Dudley and Norvaiša (1998, 1999), Mikosch and Norvaiša (2000), Lyons (1998) and so on.

For a $p \in (0, \infty)$ and a function $f$ defined on the interval $[0, T]$ and taking values in $(X, \rho)$, its *p*-variation is $v_p(f) := \sup\{\sum_{k=0}^{m-1} \rho(f(t_{k+1}), f(t_k))^p : 0 =$ $t_0 < t_1 < \cdots < t_m = T, m = 1, 2, \ldots\}$. For many important processes, conditions for the almost sure boundedness of $p$-variation of sample paths on bounded intervals are already known: Lévy (1940) showed that, for Brownian motion, the $p$-variation is bounded iff $p > 2$, a result which was refined by Taylor (1972). Later, Monroe (1978) proved that every semimartingale is equivalent to a time change of a Brownian motion. In particular, this implies that continuous-time martingales also have almost surely bounded $p$-variation for any $p > 2$. In the case of martingales with nonconstant sample paths, this result is sharp [e.g., Dudley and Norvaiša (1998), Theorem 5.6].

The $p$-variation of Gaussian processes was investigated by Jain and Monrad (1983). The same problem for Lévy processes was addressed in the papers by Fristedt and Taylor (1973) and Bretagnolle (1972). Jacob and Schilling (2001) considered more general Feller processes generated by pseudodifferential operators with uniformly bounded coefficients.

Feller processes and, in particular, Lévy processes are Markov but the methods used to establish the boundedness of the $p$-variation of their paths differ from paper to paper and do not extend to general strong Markov processes; for example, Jacob and Schilling [(2001), Theorem 2.13] require a strong assumption $p \leq 1$.

We, on the other hand, rely on the properties of general strong Markov processes and their transition probabilities. The following definition will be used throughout.

DEFINITION 1.1. Let $\beta \geq 1$ and $\gamma > 0$. We say that a Markov process $\xi_t, t \in [0, T]$, belongs to the class $\mathcal{M}(\beta, \gamma)$ ($\mathcal{M}$ for Markov) if there exist constants $a_0 > 0$ and $K > 0$ such that, for all $h \in [0, T]$ and $a \in (0, a_0]$,

$$\alpha(h, a) \leq K \frac{h^\beta}{a^\gamma}. \tag{1.1}$$

REMARK 1.2. The function $\alpha(h, a)$ is nonincreasing in $a$ for each $h$. Thus, a process $\xi_t, t \in [0, T]$, belongs to the class $\mathcal{M}(\beta, \gamma)$ iff there exist constants $a_0 > 0$ and $K > 0$ such that, for all $h \in [0, T]$ and $a > 0$,

$$\alpha(h, a) \leq K \frac{h^\beta}{(a \wedge a_0)^\gamma}, \tag{1.2}$$

because if (1.2) holds for $a = a_0$, it will also hold for any $a > a_0$.

Note that condition (1.2) and Kinney's result allow us to choose a version of $\xi_t$ with càdlàg paths. We will use this version throughout. Moreover, there exists a random variable $M(\omega)$ such that

$$\sup_{t \in [0,T]} \rho(\xi_0(\omega), \xi_t(\omega)) \leq M(\omega) < \infty \qquad \text{a.s.}$$

Our main result is

THEOREM 1.3. *Let $\xi_t, t \in [0, T]$, be a strong Markov process with values in a complete separable metric space $(X, \rho)$. Suppose $\xi_t$ belongs to the class $\mathcal{M}(\beta, \gamma)$. Then for any $p > \gamma/\beta$, the p-variation $v_p(\xi)$ of $\xi_t$ is finite almost surely.*

REMARK 1.4. For symmetric stable Lévy processes of index $\alpha \in (0, 2]$, the condition $p > \gamma/\beta$ will be shown in Section 4, at the end of the paper, to apply with $\gamma/\beta = \alpha$. Moreover, by results of Lévy (1940) and Blumenthal and Getoor (1960), the paths of such processes have unbounded *p*-variation for $p = \alpha$. Hence, the condition on $p$ in Theorem 1.3 is sharp.

As condition (1.2) indicates, we will concentrate on the "small" oscillations of the trajectories. Let

$$\nu_b(\omega) := \sup\{k : \exists \{t_i\}_{i=1}^{2k}, t_1 < t_2 \leq t_3 < \cdots < t_{2k}, \rho(\xi_{t_{2j}}, \xi_{t_{2j-1}})(\omega) > b, j \leq k\},$$

that is, $\nu_b$ is the random number of oscillations of size $> b$ of $\xi_t$ over nonoverlapping intervals.

REMARK 1.5. The paths of $\xi_t$ are chosen to be càdlàg, so the supremum in the definition of $\nu_b$ can be taken over $t_i$'s from some countable dense subset $\Gamma \subset [0, T]$; thus $\nu_b$ is a measurable random variable. Moreover, $\nu_b < \infty$ almost surely. Hence the contribution to $v_p(\xi(\omega))$ of oscillations $> b$ is at most $((2M)^p \nu_b)(\omega) < \infty$ with probability 1; that is, "large" oscillations can be handled easily.

The abundance of "small" oscillations will be handled via stopping-time techniques in Section 2.

Later we will need an Ottaviani-type inequality. Recall the following result from Gikhman and Skorohod [(1974), page 420]. Similar inequalities can be found in Blumenthal [(1957), Lemmas 2.1 and 2.2].

LEMMA 1.6. *Let $\xi_t$, $t \in [0, T]$, be a separable Markov process taking values in a complete separable metric space $(X, \rho)$. Then for any $h > 0$ and $M > 0$ such that $\alpha(h, M/2) < 1$, and for any $t \in [0, T]$,*

$$(1.3) \qquad P\left(\sup_{s \in [t, (t+h) \wedge T]} \rho(\xi_t, \xi_s) > M\right) \leq \frac{P(\rho(\xi_t, \xi_{(t+h) \wedge T}) > M/2)}{1 - \alpha(h, M/2)}.$$

**2. Markov times and related results.** Consider the natural filtration $\mathcal{F}^\xi$ generated by the process $\xi$, that is, $\mathcal{F}^\xi = \{\mathcal{F}^\xi_t, t \in [0, T]\}$, where $\mathcal{F}^\xi_t = \sigma(\xi_u, 0 \leq u \leq t) \subset \mathcal{F}$. Recall that a random variable $\tau : \Omega \to [0, \infty]$ is an $\mathcal{F}^\xi$-Markov time iff, for all $u \in [0, T]$, $\{\tau < u\} \in \mathcal{F}^\xi_u$. If $\tau$ is an $\mathcal{F}^\xi$-Markov time, define $\mathcal{F}_\tau := \{A : A \cap \{\tau < u\} \in \mathcal{F}_u^\xi, u \in [0,T]\}$. Also set $\mathcal{F}_0 := \{\varnothing, \Omega\}$. Furthermore, for any $0 \leq a < b \leq T$, let

$$R(a,b) := \sup_{a \leq s \leq t \leq b} \rho(\xi_s, \xi_t) = \sup_{s,t \in (\mathbb{Q} \cap [a,b]) \cup \{b\}} \rho(\xi_s, \xi_t),$$

since $\xi_t$ has càdlàg paths. Hence $R(a,b)$ is $\mathcal{F}_b^\xi$-measurable. Moreover, for any sequence $0 < a_n \downarrow a \leq b \leq T$, we have $R(a_n, b) \uparrow R(a,b)$ as $n \to \infty$ since the intervals $[a_n, b]$ are expanding and $\xi_t$ is right-continuous.

For any $r = 0, \pm 1, \pm 2, \ldots$, define $M_r := 2^{-r-1}$ and let $\{\tau_{l,r}\}$, $l = 0, 1, 2, \ldots$, be the sequence of random times defined as follows:

$$\tau_{0,r} := 0, \qquad \tau_{l,r} := \begin{cases} \inf\{t \in [\tau_{l-1,r}, T] : R(\tau_{l-1,r}, t) > M_r\}, \\ T+1, \quad \text{if the set above is empty.} \end{cases}$$

Each $\{\tau_{l,r}\}$ is an $\mathcal{F}^\xi$-Markov time. This is obvious for $l = 0$. Suppose that we have shown that $\{\tau_{l,r}\}$ are $\mathcal{F}^\xi$-Markov times for all $l = 0, 1, \ldots, k$. Now assume for a moment that $\tau_{k,r} < s \leq T$ and we have a sequence $\{a_n\}_{n=1}^\infty$ such that $a_n \downarrow \tau_{k,r}$ as $n \to \infty$. Then $R(a_n, s) \uparrow R(\tau_{k,r}, s)$ as $n \to \infty$ and so

$$\{R(\tau_{k,r}, s) > M_r\} = \bigcup_{n=1}^\infty \bigcap_{m=n}^\infty \{R(a_m, s) > M_r\}.$$

This observation helps to justify the second equality in the following: for any $t \leq T$, we have

$$\{\tau_{k+1,r} < t\} = \bigcup_{s < t, s \in \mathbb{Q}} \left(\{\tau_{k,r} < s\} \cap \{R(\tau_{k,r}, s) > M_r\}\right)$$

$$= \bigcup_{s < t, s \in \mathbb{Q}} \bigcup_{n=1}^\infty \bigcap_{m=n}^\infty \bigcup_{l=0}^{2^m - 1} \left(\{\tau_{k,r} \in [a_{l-1,m}, a_{l,m})\} \cap \{R(a_{l,m}, s) > M_r\}\right),$$

where $a_{l,m} := s(l+1)2^{-m}$ and the aforementioned $a_n$ is defined by

$$a_n := \begin{cases} a_{l,n}, & \text{if } \tau_{k,r} \in [a_{l-1,n}, a_{l,n}) \text{ for some } l = 0, 1, \ldots, 2^n - 1, \\ T+1, & \text{if } \tau_{k,r} \geq s. \end{cases}$$

Moreover, $\{\tau_{k,r} \in [a_{l-1,m}, a_{l,m})\} \in \mathcal{F}_{a_{l,m}}^\xi \subset \mathcal{F}_s^\xi$ since $\tau_{k,r}$ is an $\mathcal{F}^\xi$-Markov time by assumption and $\{R(a_{l,m}, s) > M_r\} \in \mathcal{F}_s^\xi$ since $a_{l,m} \leq s$ for all $l = 1, 2, \ldots, 2^m - 1$ and $m = 1, 2, \ldots$. Thus $\{\tau_{k+1,r} < t\} \in \mathcal{F}_t^\xi$; that is, $\tau_{k+1,r}$ is an $\mathcal{F}^\xi$-Markov time.

Let $\zeta_{i,r} := \tau_{i,r} - \tau_{i-1,r}$ for $i = 1, 2, \ldots$ and let $i_0$ be the smallest integer $i$ such that $\tau_{i,r} \geq T$ if at least one such $i$ exists; otherwise set $i_0 = +\infty$. Since $\xi_t$ is taken to be right-continuous, we have two cases:

*Case $i_0 < +\infty$.* Then $0 = \tau_{0,r} < \tau_{1,r} < \cdots < \tau_{i_0-1,r} < T \leq \tau_{i_0,r} \leq \tau_{i_0+1,r} = \tau_{j,r} = T+1$ for $j \geq i_0 + 1$ almost surely. Furthermore, $\zeta_{i,r} > 0$ for $i \leq i_0$, $\zeta_{i,r} = 0$ for $i > i_0 + 1$ and

$$\zeta_{i_0+1,r} = \begin{cases} 1, & \text{if } \tau_{i_0,r} = T, \\ 0, & \text{if } \tau_{i_0,r} = T+1, \end{cases}$$

almost surely.

*Case $i_0 = +\infty$.* Then $0 = \tau_{0,r} < \tau_{1,r} < \cdots < \tau_{j,r} < \cdots$ and $\zeta_{j,r} > 0$ for all $j \geq 1$ almost surely.

Before proceeding any further, let us clarify that for any random variables $X$ and $Y$ and a measurable set $A$, the expression "$X \leq Y$ *almost surely on the set $A$*" will mean $P(\{X \leq Y\} \cap A) = P(A)$. If $P(A) > 0$, this is equivalent to $P(X \leq Y|A) = 1$.

The next lemma is the first step towards bounding the average number of oscillations of a given "small" size that a trajectory can have.

LEMMA 2.1. *Let $r$ be any integer and let $u_0 \in [0, T \wedge 1]$ be such that $\alpha(u_0, M_{r+2}) < 1$. Then for any $i = 1, 2, \ldots$ and $u \in [0, u_0]$, almost surely on the set $\{\tau_{i-1,r} < T\}$,*

$$P(\zeta_{i,r} \leq u | \mathcal{F}_{\tau_{i-1,r}}) \leq \frac{\alpha(u, M_{r+2})}{1 - \alpha(u, M_{r+2})}.$$

PROOF. First of all, notice that for any $i = 1, 2, \ldots, i_0$ and $u \in [0, 1]$,

$$\begin{aligned}
&\{\zeta_{i,r} \leq u, \tau_{i-1,r} < T\} \\
&\subset \left\{ \sup_{\tau_{i-1,r} \leq t \leq (\tau_{i-1,r}+u) \wedge T} R(\tau_{i-1,r}, t) \geq M_r, \tau_{i-1,r} < T \right\} \\
(2.1) \quad &= \left\{ \sup_{s,t \in [\tau_{i-1,r}, (\tau_{i-1,r}+u) \wedge T]} \rho(\xi_s, \xi_t) \geq M_r, \tau_{i-1,r} < T \right\} \\
&\subset \left\{ \sup_{t \in [\tau_{i-1,r}, (\tau_{i-1,r}+u) \wedge T]} \rho(\xi_{\tau_{i-1,r}}, \xi_t) \geq M_{r+1}, \tau_{i-1,r} < T \right\}.
\end{aligned}$$

Take any set $A \in \mathcal{F}_{\tau_{i-1,r}}$. Then the definition of conditional expectation, the strong Markov property, (2.1) and Lemma 1.6 imply

$$\begin{aligned}
&P(\{\zeta_{i,r} \leq u\} \cap A \cap \{\tau_{i-1,r} < T\}) \\
&= \int_{A \cap \{\tau_{i-1,r} < T\}} P(\zeta_{i,r} \leq u | \mathcal{F}_{\tau_{i-1,r}}) \, dP \\
&= \int_{\{\tau_{i-1,r} < T\}} P(\zeta_{i,r} \leq u | \tau_{i-1,r}) P(A | \tau_{i-1,r}) \, dP
\end{aligned}$$

$$\leq \int_{\{\tau_{i-1,r}<T\}} P\bigg(\sup_{\tau_{i-1,r}\leq s\leq (\tau_{i-1,r}+u)\wedge T} \rho(\xi_{\tau_{i-1,r}},\xi_s)\geq M_{r+1}\bigg|\tau_{i-1,r}\bigg)$$
$$\times P(A|\tau_{i-1,r})\,dP$$
$$\leq \frac{\alpha(u,M_{r+2})}{1-\alpha(u,M_{r+2})}P(A\cap\{\tau_{i-1,r}<T\}),$$

provided $u\leq u_0$. Now, for each $\varepsilon>0$, let

$$A_\varepsilon:=A_\varepsilon(i,r,u):=\bigg\{P(\zeta_{i,r}\leq u|\mathcal{F}_{\tau_{i-1,r}})-\frac{\alpha(u,M_{r+2})}{1-\alpha(u,M_{r+2})}\geq\varepsilon\bigg\}.$$

Obviously $A_\varepsilon\in\mathcal{F}_{\tau_{i-1,r}}$, and the argument above with $A_\varepsilon$ in place of $A$ shows that

$$P(\{\zeta_{i,r}\leq u\}\cap A_\varepsilon\cap\{\tau_{i-1,r}<T\})\leq \frac{\alpha(u,M_{r+2})}{1-\alpha(u,M_{r+2})}P(A_\varepsilon\cap\{\tau_{i-1,r}<T\}).$$

On the other hand, from the definition of $A_\varepsilon$ we get

$$P(\{\zeta_{i,r}\leq u\}\cap A_\varepsilon\cap\{\tau_{i-1,r}<T\})$$
$$=\int_{A_\varepsilon\cap\{\tau_{i-1,r}<T\}}P(\zeta_{i,r}\leq u|\mathcal{F}_{\tau_{i-1,r}})\,dP$$
$$\geq\bigg(\varepsilon+\frac{\alpha(u,M_{r+2})}{1-\alpha(u,M_{r+2})}\bigg)P(A_\varepsilon\cap\{\tau_{i-1,r}<T\}).$$

Combining inequalities yields $P(A_\varepsilon\cap\{\tau_{i-1,r}<T\})=0$ for all $\varepsilon>0$. Thus, for all $\varepsilon>0$, almost surely on the set $\{\tau_{i-1,r}<T\}$,

$$P(\zeta_{i,r}\leq u|\mathcal{F}_{\tau_{i-1,r}})<\varepsilon+\frac{\alpha(u,M_{r+2})}{1-\alpha(u,M_{r+2})}.$$

Now let $\varepsilon\downarrow 0$ to get the desired inequality. □

REMARK 2.2. For any $\omega$ such that $\tau_{i-1,r}(\omega)=T$, we always have $\tau_{i,r}(\omega)=T+1$ and $\zeta_{i,r}(\omega)=1$. Hence on the set $\{\tau_{i-1,r}=T\}$,

$$P(\zeta_{i,r}\leq u|\mathcal{F}_{\tau_{i-1,r}})=\begin{cases}1,&\text{if }u\geq 1,\\ 0,&\text{otherwise.}\end{cases}$$

REMARK 2.3. Similarly, if $\omega$ is such that $\tau_{i-1,r}(\omega)=T+1$, then also $\tau_{i,r}(\omega)=T+1$ and $\zeta_{i,r}(\omega)=0$. Thus on the set $\{\tau_{i-1,r}=T+1\}$, we have, for all $u>0$,

$$P(\zeta_{i,r}\leq u|\mathcal{F}_{\tau_{i-1,r}})=1.$$

Here is a bound for the conditional expectation of the desired quantity.

LEMMA 2.4. *Let $r$ be any integer. For any $i = 1, 2, \ldots$ and $T_0 \in (0, T \wedge 1]$ such that $\alpha(T_0, M_{r+2}) \leq 1/2$, we have*

$$E(e^{-\zeta_{i,r}}|\mathcal{F}_{\tau_{i-1,r}}) \leq \begin{cases} e^{-T_0} + 2\int_0^{T_0} \alpha(x, M_{r+2})e^{-x}\,dx, & \text{a.s. on } \{\tau_{i-1,r} < T\}, \\ 1, & \text{on } \{\tau_{i-1,r} \geq T\}. \end{cases}$$

PROOF. By Remarks 2.2 and 2.3, on the set $\{\tau_{i-1,r} \geq T\}$, $\zeta_{i,r} = 0$ or 1, hence trivially $e^{-\zeta_{i,r}} \leq 1$. So consider what happens on $\{\tau_{i-1,r} < T\}$. Using Theorem 10.2.5 of Dudley (2002), integration by parts, Lemma 2.1 and the assumption $\alpha(T_0, M_{r+2}) \leq 1/2$,

$$\begin{aligned}
E(e^{-\zeta_{i,r}}|\mathcal{F}_{\tau_{i-1,r}}) &= \int_0^{T+1} e^{-x}\,dP(\zeta_{i,r} \leq x|\mathcal{F}_{\tau_{i-1,r}}) \\
&\leq e^{-(T+1)} + \int_{T_0}^{T+1} e^{-x}\,dx + \int_0^{T_0} P(\zeta_{i,r} \leq x|\mathcal{F}_{\tau_{i-1,r}})e^{-x}\,dx \\
&\leq e^{-T_0} + 2\int_0^{T_0} \alpha(x, M_{r+2})e^{-x}\,dx.
\end{aligned}$$
□

The form of condition (1.2) and the bound of Lemma 2.4 indicate that we need to deal with the incomplete gamma function which, for any $a > 0$ and $x \geq 0$, is defined as

$$\gamma(a, x) := \int_0^x u^{a-1}e^{-u}\,du.$$

Here are some needed facts about $\gamma(a, x)$.

LEMMA 2.5. (i) *For any $a > 0$ and $x \geq 0$,*

$$\gamma(a, x) = \sum_{k=0}^{\infty} \frac{(-1)^k x^{k+a}}{k!(k+a)}.$$

(ii) *If $0 \leq x < 3(3+a)(2+a)^{-1}$, then*

$$\gamma(a, x) \leq \frac{x^a}{a}\left(1 - \frac{a}{a+1}x + \frac{a}{2(a+2)}x^2\right),$$

*with strict inequality holding for $x > 0$.*

(iii) *If $a \geq 1$, then*

$$0 < a\gamma(a, 1) - \frac{1}{e} = \sum_{k=0}^{\infty} \frac{(-1)^k}{k!(k+a+1)} < \frac{7}{24}.$$

PROOF. (i) This series is well known [e.g., Davis (1970), equations (6.5.4) and (6.5.29), pages 260 and 262].

(ii) For $k \geq 2$ and $x$ as specified, the ratio of absolute values of the successive terms in the series of part (i) is

$$\frac{x(k+a)}{(k+1)(k+a+1)} < 1.$$

Since the series alternates and the third term is negative, we can discard all terms starting with the third to get the needed inequality.

(iii) From the absolute convergence of the series for $\gamma(a,x)$ and for $e^{-x}$,

$$a\gamma(a,1) - \frac{1}{e} = \sum_{k=0}^{\infty} \frac{(-1)^k a}{k!(k+a)} - \sum_{k=0}^{\infty} \frac{(-1)^k}{k!}$$

$$= \sum_{k=0}^{\infty} \frac{(-1)^k}{k!} \frac{-k}{k+a} = \sum_{k=0}^{\infty} \frac{(-1)^k}{k!(k+a+1)}.$$

Since the series alternates plus, minus, and so on, and the terms decrease in absolute value, two terms provide a lower bound, whereas three terms give an upper bound:

$$a\gamma(a,1) - \frac{1}{e} > \frac{1}{a+1} - \frac{1}{2+a} = \frac{1}{(1+a)(2+a)} > 0,$$

and for $a \geq 1$,

$$a\gamma(a,1) - \frac{1}{e} < \frac{1}{(a+1)(a+2)} + \frac{1}{2(a+3)} \leq \frac{1}{6} + \frac{1}{8} = \frac{7}{24}. \qquad \square$$

Combining the last two lemmas, we get

COROLLARY 2.6. *Assume that a strong Markov process $\xi_t, t \in [0,T]$, belongs to the class $\mathcal{M}(\beta,\gamma)$ and satisfies the conditions of Lemma 2.4. Then almost surely on $\{\tau_{i-1,r} < T\}$,*

$$E(e^{-\zeta_{i,r}}|\mathcal{F}_{\tau_{i-1,r}}) \leq \begin{cases} \beta\gamma(\beta,T_r)T_r^{-\beta}, & \text{if } T_r < 1, \\ e^{-1} + \frac{7}{24}, & \text{if } T_r = 1, \end{cases}$$

*where $T_r = \min\{((M_{r+2} \wedge a_0)^\gamma/(2K))^{1/\beta}, T, 1\}$.*

PROOF. First notice that condition (1.2) implies

$$\alpha(T_r, M_{r+2}) \leq KT_r^\beta(M_{r+2} \wedge a_0)^{-\gamma} \leq 1/2.$$

Thus, we can apply Lemma 2.4 with $T_0 = T_r$ to obtain, on $\{\tau_{i-1,r} < T\}$,

$$(2.2) \qquad E(e^{-\zeta_{i,r}}|\mathcal{F}_{\tau_{i-1,r}}) \leq e^{-T_r} + \frac{2K}{(M_{r+2} \wedge a_0)^\gamma} \int_0^{T_r} x^\beta e^{-x}\,dx.$$

Now integration by parts yields

$$0 \leq \gamma(\beta+1, T_r) = \int_0^{T_r} x^\beta e^{-x}\, dx$$
(2.3)
$$= -x^\beta e^{-x}\big|_0^{T_r} + \beta \int_0^{T_r} x^{\beta-1} e^{-x}\, dx$$
$$= -T_r^\beta e^{-T_r} + \beta \gamma(\beta, T_r),$$

and $\gamma(1, T_r) = 1 - e^{-T_r}$. Now consider three cases.

*Case* 1. $T_r^\beta = 1 \leq (M_{r+2} \wedge a_0)^\gamma/(2K)$. Then from (2.2), (2.3) and part (iii) of Lemma 2.5, we have

$$E(e^{-\zeta_{i,r}}|\mathcal{F}_{\tau_{i-1,r}}) \leq e^{-1} + \frac{2K}{(M_{r+2} \wedge a_0)^\gamma}(-e^{-1} + \beta\gamma(\beta, 1)) \leq e^{-1} + \frac{7}{24} < 0.660.$$

*Case* 2. $T_r^\beta = (M_{r+2} \wedge a_0)^\gamma/(2K)$. Then from (2.2) and (2.3),

$$E(e^{-\zeta_{i,r}}|\mathcal{F}_{\tau_{i-1,r}}) \leq e^{-T_r} + \frac{1}{T_r^\beta}(-T_r^\beta e^{-T_r} + \beta\gamma(\beta, T_r)) = \beta\gamma(\beta, T_r)T_r^{-\beta}.$$

*Case* 3. $T_r^\beta = T^\beta < \min\{1, (M_{r+2} \wedge a_0)^\gamma/(2K)\}$. Then by (2.2) and (2.3) again,

$$E(e^{-\zeta_{i,r}}|\mathcal{F}_{\tau_{i-1,r}}) \leq e^{-T} + \frac{1}{T^\beta}(-T^\beta e^{-T} + \beta\gamma(\beta, T)) = \beta\gamma(\beta, T)T^{-\beta}.$$

And so in all cases we have the stated bounds. □

Now we are ready to state and prove the following crucial lemma.

LEMMA 2.7. *Let $r$ be any integer. For any $j = 1, 2, \ldots,$*

$$P(\tau_{j,r} \leq T) \leq e^T \begin{cases} (\beta\gamma(\beta, T_r)T_r^{-\beta})^j, & \text{if } T_r < 1, \\ (e^{-1} + \frac{7}{24})^j, & \text{if } T_r = 1. \end{cases}$$

PROOF. Properties of conditional expectations, Markov's inequality and Corollary 2.6 applied $j$ times yield, for any $r$ such that $T_r < 1$,

$$P(\tau_{j,r} \leq T) = P(\tau_{j,r} \leq T, \tau_{j-1,r} < T)$$
$$= E(\mathbf{1}_{\{\tau_{j-1,r} < T\}} P(\tau_{j,r} \leq T | \mathcal{F}_{\tau_{j-1,r}}))$$
$$= E(\mathbf{1}_{\{\tau_{j-1,r} < T\}} P(e^{-\tau_{j,r}} \geq e^{-T} | \mathcal{F}_{\tau_{j-1,r}}))$$
$$\leq e^T E(\mathbf{1}_{\{\tau_{j-1,r} < T\}} E(e^{-\tau_{j,r}} | \mathcal{F}_{\tau_{j-1,r}}))$$
$$= e^T E\left(\prod_{i=1}^{j-1} (\mathbf{1}_{\{\tau_{i,r} < T\}} e^{-\zeta_{i,r}}) E(e^{-\zeta_{j,r}} | \mathcal{F}_{\tau_{j-1,r}})\right)$$

$$\leq e^T \frac{\beta\gamma(\beta, T_r)}{T_r^\beta} E\left(\prod_{i=1}^{j-2}(\mathbf{1}_{\{\tau_{i,r}<T\}} e^{-\zeta_{i,r}}) E(e^{-\zeta_{j-1,r}}|\mathcal{F}_{\tau_{j-2,r}})\right)$$

$$\leq \cdots \leq e^T \left(\frac{\beta\gamma(\beta, T_r)}{T_r^\beta}\right)^j.$$

The proof of the second case is analogous. □

**3. Proof of Theorem 1.3.** Let $PP := \{\kappa : \kappa = \{t_i : 0 = t_0 < t_1 < \cdots < t_{m_\kappa} = T\}\}$ be the set of all point partitions of $[0, T]$. For $p \geq 1$, $f : [0, T] \to X$, where as before $(X, \rho)$ is a complete separable metric space, and for $\kappa \in PP$, $\kappa = \{t_i\}_{i=0}^{m_\kappa}$, let $f_{i,\kappa} := \rho(f(t_i), f(t_{i-1})), i = 1, \ldots, m_\kappa$, and $s_p(f, \kappa) := \sum_{i=1}^{m_\kappa} f_{i,\kappa}^p$. In such notation the $p$-variation of $f$ is $v_p(f) = \sup_{\kappa \in PP} s_p(f, \kappa)$.

We now classify oscillations of a strong Markov process $\xi_t, t \in [0, T]$, which belongs to the class $\mathcal{M}(\beta, \gamma)$. For any integer $r$, recall $M_r = 2^{-r-1}$ and let $\kappa \in PP, \kappa = \{t_i\}_{i=0}^{m_\kappa}$ be an arbitrary point partition. Define the random sets

$$K_r(\omega) := K_r(\omega, \kappa) := \{k : 1 \leq k \leq m_\kappa, M_r \leq \xi_{k,\kappa} < M_{r-1}\}.$$

Let $r_1$ be the largest integer less than or equal to $-(\log_2 a_0 + 3)$, so that $M_{r+2} \geq a_0$ for all $r \leq r_1$. Let $B$ be the set of such $\omega \in \Omega$ that $\xi_t(\omega)$ is càdlàg. Then our assumptions imply $P(B) = 1$.

Now recall $\nu_b$ as defined before Remark 1.5. For any $\omega \in B$, let $\nu_0(\omega) := \nu_{a_0/2}(\omega)$, so that $\nu_0$ counts all oscillations $> a_0/2$, including all those $\geq a_0$. Then write

$$v_p(\xi(\omega)) \leq \sum_{r>r_1} \sup_{\kappa \in PP} \sum_{k \in K_r(\omega)} \xi_{k,\kappa}^p(\omega) + \sup_{\kappa \in PP} \sum_{r \leq r_1} \sum_{k \in K_r(\omega)} \xi_{k,\kappa}^p(\omega)$$

(3.1)
$$\leq \sum_{r>r_1} 2^{-rp} \sup_{\kappa \in PP} \sum_{k \in K_r(\omega)} 1 + ((2M)^p \nu_0)(\omega)$$

$$=: S_1 + ((2M)^p \nu_0)(\omega).$$

To establish that $v_p(\xi)$ is bounded almost surely for $p > \gamma/\beta$, it is sufficient to bound

$$E \sup_{\kappa \in PP} \sum_{k \in K_r(\omega)} 1 =: EY_r$$

for $r > r_1$, and to show that $S_1$ converges for $p > \gamma/\beta$.

LEMMA 3.1. *Let $r$ be any integer and let $\alpha(x, a)$ and $T_r$ be as in Corollary 2.6. Suppose that $r > r_1$. Then*

$$EY_r \leq \begin{cases} 4T_r^{-1} e^T, & \text{if } T_r < 1, \\ 1.95 e^T, & \text{if } T_r = 1. \end{cases}$$

PROOF. First write $EY_r = \sum_{j=1}^{\infty} P(Y_r \geq j)$. Since $r > r_1$, we have $M_{r+2} < a_0$ and $T_r = \min\{(M_{r+2}^\gamma/(2K))^{1/\beta}, T, 1\} < 1$. So by Lemma 2.7, we get

$$P(Y_r \geq j) \leq P(\tau_{j,r} \leq T) \leq e^T \left\{ \frac{\beta\gamma(\beta, T_r)}{T_r^\beta} \right\}^j.$$

Now Lemma 2.5(ii) yields $\beta\gamma(\beta, T_r) < T_r^\beta$ and, since $T_r < 1$, we also have

$$\frac{1}{T_r^\beta - \beta\gamma(\beta, T_r)} \leq \frac{1}{\beta/(\beta+1)T_r^{\beta+1} - \beta/(2(\beta+2))T_r^{\beta+2}} \leq \frac{2(\beta+1)}{\beta T_r^{\beta+1}} \leq \frac{4}{T_r^{\beta+1}}.$$

Thus for $T_r < 1$,

$$EY_r \leq e^T \sum_{j=1}^{\infty} \left\{ \frac{\beta\gamma(\beta, T_r)}{T_r^\beta} \right\}^j$$

$$= e^T \frac{\beta\gamma(\beta, T_r)}{T_r^\beta} \left\{ 1 - \frac{\beta\gamma(\beta, T_r)}{T_r^\beta} \right\}^{-1}$$

$$\leq 4e^T \frac{\beta\gamma(\beta, T_r)}{T_r^{\beta+1}} \leq \frac{4e^T}{T_r}.$$

If $T_r = 1$, then Lemma 2.7 gives, for each $j$,

$$P(Y_r \geq j) \leq P(\tau_{j,r} \leq T) \leq e^T (e^{-1} + \tfrac{7}{24})^j < e^T (0.66)^j.$$

Thus in this case

$$EY_r \leq e^T \sum_{j=1}^{\infty} (0.66)^j = e^T \tfrac{0.66}{0.34} < 1.95 e^T. \qquad \square$$

Returning to the proof of Theorem 1.3, recall the definition of $r_1$. It satisfies the following inequalities:

(3.2) $$r_1 \leq -(\log_2 a_0 + 3) < r_1 + 1.$$

Now choose any $\varepsilon > 0$ and let $N > 0$. Then from (3.1) we get

$$P(v_p(\xi) > N) = P(\{v_p(\xi) > N\} \cap B) \leq P_1 + P_2 + P_3,$$

where, using Markov's inequality and Lemma 3.1,

$$P_1 := P\left(S_1 > \frac{N}{2}\right) \leq \frac{2}{N} \sum_{r > r_1} 2^{-rp} EY_r$$

$$\leq \frac{2e^T}{N} \sum_{r > r_1} 2^{-rp} \{4T_r^{-1} \mathbf{1}_{\{T_r < 1\}} + (1.95) \mathbf{1}_{\{T_r = 1\}}\}.$$

If $r > r_1$, then $M_{r+2} = 2^{-r-3} < a_0$ and $T_r = \min\{(M_{r+2}^\gamma/(2K))^{1/\beta}, T, 1\}$. And so

$$T_r^{-1} = \left(\min\left\{\left(\frac{M_{r+2}^\gamma}{2K}\right)^{1/\beta}, T, 1\right\}\right)^{-1}$$

$$= \max\left\{1, \frac{1}{T}, \frac{(2K)^{1/\beta}}{M_{r+2}^{\gamma/\beta}}\right\}$$

$$\leq 1 + \frac{1}{T} + (K2^{1+\gamma(r+3)})^{1/\beta}.$$

Hence

$$P_1 \leq \frac{2e^T}{N}\left\{\left(6+\frac{4}{T}\right)\sum_{r=r_1+1}^\infty 2^{-rp} + K^{1/\beta}2^{2+(3\gamma+1)/\beta}\sum_{r=r_1+1}^\infty 2^{-r(p-\gamma/\beta)}\right\}$$

$$= \frac{2e^T}{N}\left\{\left(6+\frac{4}{T}\right)\frac{2^{-(r_1+1)p}}{(1-2^{-p})} + \frac{K^{1/\beta}2^{2+(3\gamma+1)/\beta-(r_1+1)(p-\gamma/\beta)}}{1-2^{-(p-\gamma/\beta)}}\right\}$$

$$=: \frac{C_1(K,T,p,\beta,\gamma,a_0)}{N} \leq \frac{\varepsilon}{3},$$

for all $N \geq N_1 := [3C_1(K,T,p,\beta,\gamma,a_0)\varepsilon^{-1}] + 1$ and $p > \gamma/\beta$. Moreover, for $N \geq N_2(\varepsilon)$ large enough,

$$P_2 := P\left((2M(\omega))^p > \sqrt{\frac{N}{2}}\right) \leq \frac{\varepsilon}{3}.$$

Similarly, for $N \geq N_3(\varepsilon)$ large enough,

$$P_3 := P\left(\nu_0(\omega) > \sqrt{\frac{N}{2}}\right) \leq \frac{\varepsilon}{3}.$$

Combining the obtained bounds and recalling Remark 1.5, for all $N \geq \max\{N_1, N_2, N_3\}$ and $p > \gamma/\beta$, we have $P(v_p(\xi) > N) \leq \varepsilon$. This implies that, for $p > \gamma/\beta$, $v_p(\xi)$ is bounded almost surely. This completes the proof of the main theorem.

**4. Examples.** In this section we show how Theorem 1.3 can be applied to real-valued symmetric $\alpha$-stable Lévy processes.

Let $X_1$ be a real-valued symmetric $\alpha$-stable random variable with a characteristic function $\phi(t) = e^{-c|t|^\alpha}$ for $t \in \mathbb{R}$, $\alpha \in (0,2]$ and some constant $c > 0$. Consider the temporally homogeneous symmetric $\alpha$-stable Lévy motion $\{X_t, t \geq 0\}$ started at $x = 0$ and with increments having distribution:

$$X_t - X_s \sim |t-s|^{1/\alpha}X_1, \qquad t,s > 0.$$

[For more on Lévy processes and their properties, see, e.g., Sato (2000).]

It is well known [e.g., Feller (1971), page 448] that, for some constant $b \in (0, \infty)$,

$$x^\alpha P(X_1 > x) \to b \qquad \text{as } x \to +\infty.$$

Thus there exists an $x_0 \in (0, \infty)$ such that

$$P(X_1 > x) \leq 2bx^{-\alpha} \qquad \text{for all } x \geq x_0.$$

For $x \in (0, x_0)$, we trivially have $P(X_1 > x) \leq x_0^\alpha x^{-\alpha}$. By assumption, $X_1$ is symmetric, and so combining the above bounds, we have

$$P(|X_1| > x) \leq 2\max\{x_0^\alpha, 2b\}x^{-\alpha} =: Kx^{-\alpha}.$$

Therefore for $h > 0$ and $a > 0$,

$$(4.1) \qquad P(|X_{t+h} - X_t| > a) = P(|X_1| > ah^{-1/\alpha}) \leq K\frac{h}{a^\alpha}.$$

Denote $A_a(x) = \{y \in \mathbb{R} : |x - y| \geq a\}$ for $a > 0$. From the spatial homogeneity of transition probability functions of Lévy processes [e.g., Sato (2000), Theorem 10.5] and (4.1), we get, for any $T > 0$, $h \in (0, T]$ and $a > 0$,

$$(4.2) \qquad \begin{aligned} \alpha(h, a) &= \sup\{P_{s,t}(x, A_a(x)) : x \in \mathbb{R}, 0 \leq s \leq t \leq (s+h) \wedge T\} \\ &= \sup\{P_{s,t}(0, A_a(0)) : 0 \leq s \leq t \leq (s+h) \wedge T\} \\ &= \sup\{P(|X_t - X_s| \geq a) : 0 \leq s \leq t \leq (s+h) \wedge T\} \\ &\leq \sup\left\{\frac{Kv}{a^\alpha} : 0 \leq v \leq h \leq T\right\} = K\frac{h}{a^\alpha}. \end{aligned}$$

Now Theorem 1.3 implies that $\alpha$-stable Lévy motion has bounded $p$-variation on every interval $[0, T]$ for $p > \alpha$. This reestablishes the "positive" parts of Theorem 9 of Lévy (1940) for $\alpha = 2$ and Theorem 4.1 of Blumenthal and Getoor (1960) for $\alpha \in (0, 2)$. Moreover, this example shows that the condition "$p > \gamma/\beta$" of Theorem 1.3, in general, cannot be replaced by $p \geq \gamma/\beta$: by the aforementioned theorems, for $p = \gamma/\beta = \alpha$, the $p$-variation of an $\alpha$-stable Lévy motion is infinite almost surely for all intervals $[0, T]$.

**Acknowledgments.** This paper contains some of the results from my Ph.D. dissertation written under the guidance of Professor Richard M. Dudley, to whom I am deeply grateful for many helpful suggestions. I would also like to thank Professors Marjorie G. Hahn and Daniel W. Stroock as well as an anonymous referee who carefully read this paper and helped to improve it.

Department of Mathematics
University of Connecticut
196 Auditorium Road, Unit 9
Storrs, Connecticut 06269
USA
e-mail: martynas@math.uconn.edu